\documentclass[11pt]{amsart}
\usepackage[T1]{fontenc}
\usepackage{lmodern}
\usepackage{amsthm,amssymb,mathtools}
\usepackage[colorlinks]{hyperref}
\usepackage[margin=1in]{geometry}
\usepackage{enumitem}
\usepackage{microtype}
\usepackage{tikz}
\usepackage{booktabs}
\usepackage{needspace}

\allowdisplaybreaks[2]
\setlist[enumerate]{label=\textup{(\roman*)},leftmargin=2em,itemsep=0.2em}
\newtheoremstyle{citedplain}{\topsep}{\topsep}{\itshape}{}{\bfseries}{.}{.5em}
  {\thmname{#1}\thmnumber{ #2}\thmnote{ \textnormal{[#3]}}}
\newcommand{\bibsource}[2]{\nocite{#1}\hyperlink{cite.#1}{#2}}
\theoremstyle{citedplain}
\newtheorem{theorem}{Theorem}[section]
\newtheorem{thmx}{Theorem}

\newtheorem{lemma}[theorem]{Lemma}
\newtheorem{proposition}[theorem]{Proposition}
\newtheorem{corollary}[theorem]{Corollary}
\newtheorem{question}[theorem]{Question}
\theoremstyle{definition}

\newtheorem{example}[theorem]{Example}
\newtheorem{remark}[theorem]{Remark}
\numberwithin{equation}{section}
\newcommand{\Z}{\mathbb Z}
\newcommand{\R}{\mathbb R}
\newcommand{\C}{\mathbb C}
\newcommand{\vs}{\vspace}

\newcommand{\cF}{\mathcal F}
\newcommand{\hh}{\widetilde \Phi}
\DeclareMathOperator{\ind}{ind}
\DeclareMathOperator{\sgn}{sgn}
\DeclareMathOperator{\supp}{supp}
\hypersetup{pdftitle={Unimodality of Betti numbers for Hamiltonian circle actions in dimension ten},
 pdfauthor={Yunhyung Cho},
 pdfsubject={Equivariant cohomology, Hamiltonian circle actions, and unimodality},
 linkcolor=red!50!black,citecolor=blue!70!black,urlcolor=blue!55!black}
\begin{document}
\title[Unimodality for Hamiltonian circle actions in dimension ten]{Unimodality of Betti numbers for Hamiltonian circle actions in dimension ten}
\author{Yunhyung Cho}
\address{Department of Mathematics Education, Sungkyunkwan University, Seoul, Republic of Korea}
\email{yunhyung@skku.edu}

\begin{abstract}
Kirwan's injectivity theorem states that the equivariant cohomology of a closed connected Hamiltonian $S^1$-manifold $(M^{2n},\omega)$ injects into that of its fixed point set. 
In this paper, we refine this theorem degree by degree: for nonnegative integers $a$ and $d$ with $2a+\max\{1,2d\}<n$, we construct a collection $\mathcal C_{a,d}$ of fixed components such that restriction to $\mathcal C_{a,d}$ is injective in degree $2a$. 
When all fixed points are isolated, this yields inequalities among the Betti numbers of $M$ in arbitrary dimension. 
In real dimension ten, we prove unimodality of the even Betti numbers, answering Tolman's question affirmatively for actions with isolated fixed points.
\end{abstract}

\maketitle

\section{Introduction}\label{sec:introduction}

Let $(M^{2n},\omega)$ be a closed connected symplectic manifold with an effective Hamiltonian circle action and moment map $\Phi$.
Kirwan's injectivity theorem \cite{Kirwan} states that an equivariant cohomology class on $M$ is determined by its restrictions to the fixed components. 
Injectivity fails if even one of the fixed components is omitted. 
In fact, the equivariant Poincar\'e dual of a fixed component $F$ has nonzero restriction to $F$ and zero restriction to every other component. 

Now, fix an even degree $2a$. We say that a collection $\mathcal C$ of fixed components \emph{detects} $H^{2a}_{S^1}(M;\R)$ if the restriction map
\[
	H^{2a}_{S^1}(M;\R)\longrightarrow \bigoplus_{F\in\mathcal C}H^{2a}_{S^1}(F;\R)
\]
is injective. 
The Morse--Bott filtration of the moment map gives one such collection, consisting of the fixed components of Morse index at most $2a$ (Proposition~\ref{prop:elementary}).
In this paper, we construct other collections indexed by pairs $(a,d)\in(\Z_{\geq 0})^2$ with $2a+\max\{1,2d\}<n$ and use them to derive new inequalities among the Betti numbers of $M$.

To state our main theorem, we use the following notation.
Denote by $\cF$ the set of connected components of $M^{S^1}$, by $\cF_+$ the subset of positive-dimensional components, and by $\mathcal S$ the set of isolated fixed points.
For each $F\in\cF$, let $i_F:F\hookrightarrow M$ be the inclusion and write 
\begin{equation}\label{eq:indices}
	f(F)=\tfrac12\dim_{\R}F,\qquad \mu(F)=\tfrac12\ind_\Phi(F),\qquad c(F)=\tfrac12\ind_{-\Phi}(F)=n-f(F)-\mu(F).
\end{equation}
For an integer $d\geq 0$, set
\begin{equation}\label{eq:sets}
	L_d=\{q\in\mathcal S:c(q)>d,\ c(q)\equiv d\pmod2\},\qquad \mathcal C_{a,d}=\cF_+\cup L_d .
\end{equation}
In other words, $\mathcal C_{a,d}$ consists of all positive-dimensional fixed components and the isolated fixed points with half-coindex $d+2,d+4,\ldots$. 
Our first main result is that $\mathcal C_{a,d}$ detects $H_{S^1}^{2a}(M)$. 
\Needspace{5\baselineskip}

\begin{thmx}\label{thm:injective}
	Let $(M^{2n},\omega)$ be a closed connected Hamiltonian $S^1$-manifold with a moment map $\Phi$.
	If $a,d\ge0$ are integers satisfying
\begin{equation}\label{eq:range}
	2a+\max\{1,2d\}<n,
\end{equation}
then restriction to $\mathcal C_{a,d}$ gives an injective map
\begin{equation}\label{eq:mixed-injection}
	H^{2a}_{S^1}(M)\longrightarrow \Bigl(\bigoplus_{F\in\cF_+}H^{2a}_{S^1}(F)\Bigr)\oplus\R^{L_d},\qquad \alpha\longmapsto \bigl((i_F^*\alpha)_{F\in\cF_+},(\alpha_q)_{q\in L_d}\bigr),
\end{equation}
where $H^*_{S^1}(q;\R)=\R[u]$, $\deg u=2$, and $i_q^*\alpha=\alpha_q u^a$.
\end{thmx}

By equivariant formality, $H^*_{S^1}(M;\R)$ is isomorphic to $H^*(M;\R)\otimes\R[u]$ as a graded $\R[u]$-module. 
This gives
\[
	\dim H^{2a}_{S^1}(M;\R)=b_0(M)+b_2(M)+\cdots+b_{2a}(M).
\]
If all fixed points are isolated, then $b_{2j}(M)$ equals the number of fixed points of half-index $j$, and hence
\[
	\dim H^{2a}_{S^1}(M;\R)=|A_a|,\qquad A_a=\{q\in\mathcal S:\mu(q)\le a\}.
\]
In this case, Theorem~\ref{thm:injective} gives $|A_a|\le|L_d|$. 
Since both sides are sums of Betti numbers, this yields inequalities among the Betti numbers of $M$, some of which are new (Section~\ref{sec:numerical}). 
We note that this inequality cannot be improved in general: for suitable $(a,d)$ it is an equality for linear circle actions on complex projective spaces (Example~\ref{ex:CPn}).

The next theorem shows that the inequality $|A_a|\le |L_d|$ can be refined at every level of the moment map when all fixed points are isolated. \vs{0.2cm}

\begin{thmx}\label{thm:levels}
Let $(M^{2n},\omega)$ be a closed connected Hamiltonian $S^1$-manifold with isolated fixed points and a moment map $\Phi$. 
For integers $a,d\ge0$ with $2a+\max\{1,2d\}<n$ and every $t\in\R$, we have 
\begin{equation}\label{eq:tail}
	\#\{p\in A_a:\Phi(p)>t\} \le \#\{q\in L_d:\Phi(q)>t\}.
\end{equation}
After removing the points common to $A_a$ and $L_d$, the following stronger inequality holds:
\begin{equation}\label{eq:strict-tail}
	\#\{p\in A_a\setminus L_d:\Phi(p)\ge t\} \le \#\{q\in L_d\setminus A_a:\Phi(q)>t\}.
\end{equation}
\end{thmx}

Recall that a sequence $x_0,x_1,\ldots,x_N$ of integers is \emph{unimodal} if there is an index $k$ such that
\[
    x_0\le\cdots\le x_k\ge\cdots\ge x_N.
\]
Unimodality arises naturally in K\"ahler geometry. 
For a compact K\"ahler manifold, the hard Lefschetz theorem \cite{Huybrechts} implies that the even and odd Betti numbers each form a unimodal sequence. 
In dimension four, a Hamiltonian circle action with isolated fixed points extends to a toric action, and the manifold admits an invariant K\"ahler structure \cite{Karshon}. 
In dimension six, Tolman \cite{TolmanNK} constructed a Hamiltonian $T^2$-manifold with isolated fixed points that admits no $T^2$-invariant K\"ahler structure.
Nevertheless, Goertsches, Konstantis and Zoller later showed that the underlying manifold admits a K\"ahler structure, which is necessarily not $T^2$-invariant \cite{GKZ20}; see also \cite{GKZ23}. 
This raises the question of whether a Hamiltonian circle action with isolated fixed points forces the underlying manifold to admit a K\"ahler structure.

At the 2005 BIRS workshop ``Moment maps in various geometries'', Karshon asked whether a closed symplectic manifold with a Hamiltonian circle action and isolated fixed points satisfies the hard Lefschetz property \cite[Problem~4.2]{BIRS}.
Tolman posed a weaker question about unimodality \cite[Problem~4.3]{BIRS}, which we state here for isolated fixed points. \vs{0.2cm}

\begin{question}[\bibsource{BIRS}{Tolman}]\label{q:tolman}
	If a closed connected symplectic manifold $(M^{2n},\omega)$ admits a Hamiltonian circle action with isolated fixed points, must its even Betti numbers satisfy
	\[
		b_0(M)\le b_2(M)\le\cdots\le b_{2\lfloor n/2\rfloor}(M)?
	\]
\end{question}

In real dimension at most six, unimodality follows from Poincar\'e duality and $b_0=1\le b_2$, the latter inequality being a consequence of $[\omega]\ne 0$.
The author and Kim \cite{CK} proved unimodality in the first nontrivial case, real dimension eight. 
Luo \cite[Corollary~5.2]{Luo} also established $b_2\le b_4$ in real dimensions eight and ten for Hamiltonian GKM manifolds whose moment maps are in general position.

For actions with isolated fixed points, unimodality is also known when the action is semifree \cite{TW}, when the moment map is index-increasing \cite{ChoIndex}, or when the moment map is constant on the fixed points of each Morse index \cite{CKSelf}.
Without these additional hypotheses, the author \cite[Theorem~2]{Cho21} proved the inequality
\begin{equation}\label{eq:alter}
    b_2+b_6+\cdots+b_{4k-2} \le b_4+b_8+\cdots+b_{4k}
\end{equation}
in real dimensions $8k$ and $8k+4$.
In particular, this gives $b_2\le b_4$ in real dimension twelve.
Unimodality remains open in arbitrary dimension.

Applying Theorem~\ref{thm:injective}, we prove unimodality in real dimension ten. 
As further consequences, we also recover unimodality in real dimension eight \cite{CK} and inequality~\eqref{eq:alter}.\vs{0.2cm}

\begin{thmx}\label{thm:ten}
	Let $(M,\omega)$ be a closed connected symplectic manifold of real dimension eight or ten. If $M$ admits a Hamiltonian circle action with isolated fixed points, then its even Betti numbers are unimodal.
\end{thmx}

\begin{proof}
Take $(a,d)=(1,0)$ when $n=4$ and $(a,d)=(1,1)$ when $n=5$.
Both pairs satisfy \eqref{eq:range}. 
In either case $A_1$ consists of the fixed points of half-index $0$ or $1$, and $L_d$ consists of those of half-index $0$ or $2$. 
Thus $|A_1|\le |L_d|$ gives $b_0(M)+b_2(M)\le b_0(M)+b_4(M)$. 
Since $M$ is connected and $[\omega]\ne0$, we also have $b_0(M)=1\le b_2(M)$. 
Poincar\'e duality now gives unimodality in both dimensions.
\end{proof}

For actions with isolated fixed points in higher dimensions, different choices of $(a,d)$ give further inequalities.
For example,
\[
    b_2\le b_{2m}\quad\text{in dimension }4m+2\ (m\ge1),\qquad
    1+b_4\le 2b_6\quad\text{in dimension sixteen}.
\]
For even $n\ge4$, the choice $(a,d)=(n/2-1,0)$ recovers \eqref{eq:alter} by Poincar\'e duality (Corollary~\ref{cor:cho}). 

Related results have also been obtained for actions with positive-dimensional fixed components. 
Lindsay \cite{Lindsay25} proved $b_2\le b_4$ in real dimension eight when neither extremal fixed component has dimension four, and unimodality in real dimension twelve when $b_2=1$ and no fixed component has dimension four or eight. 
More recently, Li and Ren \cite[Section~7]{LiRen26} obtained further restrictions on the Betti numbers using the signature. 
In particular, if every fixed component is either a point or has real dimension congruent to $2$ modulo $4$, and if $b_2=1$, then $b_4\le b_6$ in real dimensions twelve and sixteen \cite[Corollary~7.8]{LiRen26}.

This paper is organized as follows. Section~\ref{sec:prelim} fixes our conventions and recalls the equivariant Morse--Bott filtration and the localization formula. 
Section~\ref{sec:restriction} proves Theorems~\ref{thm:injective} and~\ref{thm:levels}. 
Section~\ref{sec:numerical} derives the resulting inequalities among the Betti numbers in arbitrary dimension.

\section{Preliminaries}\label{sec:prelim}

In this section, we review the basic tools of equivariant cohomology for Hamiltonian circle actions: the Atiyah--Bott--Berline--Vergne (ABBV) localization formula, the equivariant Morse--Bott filtration, and the associated Morse-theoretic bases. 
Throughout the paper, all cohomology groups have real coefficients.

Let $(M^{2n},\omega)$ be a closed connected symplectic manifold equipped with a Hamiltonian circle action with a moment map $\Phi$. 
We use the sign convention $d\Phi=-\iota_X\omega$ for a moment map $\Phi$. Here, $X$ denotes the fundamental vector field of the $S^1$-action. 
Choose an $S^1$-invariant almost complex structure compatible with $\omega$. 
For each fixed component $F$, the normal bundle of $F$ splits as $N_F^-\oplus N_F^+$ according to the signs of the weights. 
With our convention, $\Phi$ decreases along the negative-weight directions, so $N_F^-$ and $N_F^+$ have complex ranks $\mu(F)$ and $c(F)$, respectively. 
Let $e_F^\pm$ denote the equivariant Euler classes of $N_F^\pm$. 
Since $S^1$ acts trivially on $F$, we have 
\[
	H^*_{S^1}(F)=H^*(F\times BS^1)=H^*(F)\otimes\R[u],
\]
where $u\in H^2(BS^1)$ is the equivariant Euler class of $\C$ with the action $\theta\cdot z=e^{2\pi i\theta}z$. 
The leading coefficient of $e_F^\pm$ in $u$ is the product of the corresponding weights, counted with multiplicity. 
Since it is a nonzero scalar, multiplication by $e_F^\pm$ on $H^*_{S^1}(F)$ is injective.

At an isolated fixed point $q$ with tangent weights $w_{q,1},\ldots,w_{q,n}$, we have
\begin{equation}\label{eq:euler}
	e_{S^1}(T_qM)=\Bigl(\prod_{\nu=1}^n w_{q,\nu}\Bigr)u^n,\qquad \lambda_q:=\Bigl(\prod_{\nu=1}^n w_{q,\nu}\Bigr)^{-1},\qquad \sgn\lambda_q=(-1)^{\mu(q)}.
\end{equation}
For $\gamma\in H^{2r}_{S^1}(M)$ and $q\in\mathcal S$, we define the \emph{coefficient} $\gamma_q\in\R$ of $\gamma$ at $q$ by
\begin{equation}\label{eq:coeff}
	i_q^*\gamma=\gamma_q u^r.
\end{equation}
We also use the equivariant symplectic class $\hh\in H^2_{S^1}(M)$, represented by $-\omega+\Phi u$ in the Cartan model with differential $d-u\iota_X$ \cite{AB}; see \cite[Chapter~9]{GS} for equivariant symplectic forms in general.
Its restriction to a fixed component is
\begin{equation}\label{eq:component-moment}
	i_F^*\hh=\Phi(F)\,u-[\omega|_F]\qquad(F\in\cF);
\end{equation}
in particular, $\hh_q=\Phi(q)$ for every $q\in\mathcal S$. \vs{0.2cm}

\begin{theorem}[\bibsource{AB}{Atiyah--Bott}, \bibsource{BV}{Berline--Vergne}]\label{thm:localization}
	Let $(M^{2n},\omega)$ be a closed connected Hamiltonian $S^1$-manifold.
	Then every $\gamma\in H^*_{S^1}(M)$ satisfies
	\begin{equation}\label{eq:localization}
		\int_M\gamma=\sum_{F\in\cF}\int_F\frac{i_F^*\gamma}{e_F^-e_F^+}.
	\end{equation}
\end{theorem}

We note that if $\gamma\in H^{2r}_{S^1}(M)$ with $0\le r<n$, then $\int_M\gamma=0$ for degree reasons. 
In addition, if such a $\gamma$ restricts to zero on every positive-dimensional fixed component, then
\[
	\sum_{q\in\mathcal S}\lambda_q\gamma_q=0.
\]

Next, we recall the Morse--Bott filtration of $\Phi$. For a critical value $t$ of $\Phi$, let
\[
	M^\pm=\{\Phi\le t\pm\epsilon\} \subseteq M
\]
where $\epsilon>0$ is small enough that $[t-\epsilon,t+\epsilon]$ contains no other critical value. 
At the minimum value $M^-=\varnothing$, and at the maximum value $M^+=M$. \vs{0.2cm}

\begin{lemma}[\bibsource{Kirwan}{Kirwan}]\label{lem:filtration}
	Let $(M^{2n},\omega)$ be a closed connected Hamiltonian $S^1$-manifold with a moment map $\Phi$.
\begin{enumerate}
\item For each critical value $t$ and every $j$, there is a short exact
sequence
\begin{equation}\label{eq:thom}
	0\longrightarrow\bigoplus_{\Phi(F)=t}H^{j-2\mu(F)}_{S^1}(F)
	\longrightarrow H^j_{S^1}(M^+)\longrightarrow H^j_{S^1}(M^-)\longrightarrow0.
\end{equation}
Here the first term is identified with $H^j_{S^1}(M^+,M^-)$ by excision and the Thom isomorphisms, the first map is the map in the long exact sequence of the pair, and the second map is restriction. 
The image of a class $\gamma$ in the $F$-summand restricts to $e_F^-\gamma$ on $F$ and to zero on the other components at level $t$. 
Moreover, for each of $M^\pm$, restriction to its fixed point set is injective \textup{(}see \cite[Lemma~3.1 and the proof of Theorem~3.2]{TW03}\textup{)}.

\item For every $F\in\cF$ and every $\gamma\in H^{j-2\mu(F)}_{S^1}(F)$, there is a class $\tau\in H^j_{S^1}(M)$ such that
\[
	i_F^*\tau=e_F^-\gamma,\qquad
	i_{F'}^*\tau=0\quad\text{if }F'\ne F\text{ and }\Phi(F')\le\Phi(F).
\]

\item For all integers $k,a\ge0$, with the convention $b_i(F)=0$ for $i<0$,
\begin{equation}\label{eq:perfection}
	b_k(M)=\sum_{F\in\cF}b_{k-2\mu(F)}(F),\qquad
	\dim H^{2a}_{S^1}(M)=\sum_{j=0}^{a}b_{2j}(M).
\end{equation}
\end{enumerate}
\end{lemma}
\vs{0.2cm}

\begin{proposition}\label{prop:elementary}
Let $(M,\omega)$ be a closed connected Hamiltonian $S^1$-manifold with moment map $\Phi$. 
For every integer $k\ge0$, restriction to the fixed components of Morse index at most $k$ gives an injection
\[
	H^k_{S^1}(M)\longrightarrow\bigoplus_{2\mu(F)\le k}H^k_{S^1}(F).
\]
\end{proposition}

\begin{proof}
Let $\alpha\in H^k_{S^1}(M)$ be nonzero. 
By Kirwan's injectivity theorem, there is a lowest critical value $t$ at which $\alpha$ restricts nontrivially to some fixed component $F$ with $\Phi(F)=t$. 
Since all restrictions of $\alpha$ below $t$ vanish, Lemma~\ref{lem:filtration}\,\textup{(i)} gives $\alpha|_{M^-}=0$.
By exactness of \eqref{eq:thom}, $i_F^*\alpha=e_F^-\gamma$ for some nonzero $\gamma\in H^{k-2\mu(F)}_{S^1}(F)$. 
Hence $2\mu(F)\le k$, so $\alpha$ restricts nontrivially to a component of Morse index at most $k$.
\end{proof}

Now suppose that all fixed points are isolated. 
Then $\Phi$ is a perfect Morse function, and \eqref{eq:perfection} gives
\begin{equation}\label{eq:morse-count}
	b_{2j}(M)=\#\{q\in\mathcal S:\mu(q)=j\},\qquad H^{2j+1}(M)=0.
\end{equation}
In this case, Lemma~\ref{lem:filtration} yields the following basis. \vs{0.2cm}

\begin{lemma}[\bibsource{GT}{Goldin--Tolman, Lemma~2.5}]\label{lem:GT}
Let $(M,\omega)$ be a closed connected Hamiltonian $S^1$-manifold with isolated fixed points and a moment map $\Phi$. 
For each fixed point $p\in\mathcal S$, there is a class $\tau_p\in H^{2\mu(p)}_{S^1}(M)$ such that
\[
	i_p^*\tau_p=e_p^-,\qquad
	i_q^*\tau_p=0\quad\text{if }q\ne p\text{ and }\Phi(q)\le\Phi(p).
\]
For every integer $a\ge0$, the classes
\[
	v_p=u^{a-\mu(p)}\tau_p,\qquad p\in A_a=\{p\in\mathcal S:\mu(p)\le a\},
\]
form a basis of $H^{2a}_{S^1}(M)$. In particular,
\begin{equation}\label{eq:dimension}
	\dim H^{2a}_{S^1}(M)=|A_a|=\sum_{j=0}^{a}b_{2j}(M).
\end{equation}
\end{lemma}

\section{Proofs of Theorems~\ref{thm:injective} and~\ref{thm:levels}}\label{sec:restriction}

Let $(M^{2n},\omega)$ be a closed connected Hamiltonian $S^1$-manifold with a moment map $\Phi$.
We use the notation of \eqref{eq:indices} and \eqref{eq:sets}. For an integer $d\ge0$, let
\[
	U_d=\{p\in\mathcal S:c(p)\le d,\ c(p)\equiv d\pmod2\}.
\]
Note that $U_d$ and $L_d$ are disjoint. \vs{0.2cm}

\begin{lemma}\label{lem:interp}
	For each $p\in U_d$, there exists a class $\eta_p\in H^{2d}_{S^1}(M)$ satisfying
\begin{equation}\label{eq:interp}
    \begin{aligned}
	    i_{q}^*\eta_p&=\delta_{pq}u^d &&(q\in U_d),\\
	    i_F^*\eta_p&=0 &&(F\ne\{p\},\ \Phi(F)\ge\Phi(p)).
    \end{aligned}
\end{equation}
\end{lemma}

\begin{proof}
	For $p\in U_d$, the negative normal bundle for $-\Phi$ at $p$ is $N_p^+$. Applying Lemma~\ref{lem:filtration}\,\textup{(ii)} to $-\Phi$ gives a class $\bar{\tau}_p\in H^{2c(p)}_{S^1}(M)$ such that
	\[
		\bar{\tau}_p|_{p} = e_p^+ = w_p^+ u^{c(p)}, \qquad \bar{\tau}_p|_F = 0 \quad (F\ne\{p\},\ \Phi(F)\ge\Phi(p)),
	\]
	where $w_p^+$ is the product of the positive weights at $p$. 
	Since $c(p)\le d$, the class $\gamma_p := \frac{1}{w_p^+}u^{d-c(p)}\bar{\tau}_p$ lies in $H^{2d}_{S^1}(M)$ and satisfies
	\[
		i_{p}^*\gamma_p =u^d, \qquad \qquad i_F^*\gamma_p =0 \quad (F\ne\{p\},\ \Phi(F)\ge\Phi(p)).
	\]
	
	Order $U_d$ by nondecreasing moment value, breaking ties arbitrarily. With $(\gamma_p)_q\in\R$ as in \eqref{eq:coeff}, define $\eta_p$ inductively by
	\[
	    \eta_p = \gamma_p-\sum_{\substack{q\in U_d\\\Phi(q)<\Phi(p)}}(\gamma_p)_q \eta_q.
	\]
	At each $q\in U_d$ with $\Phi(q)<\Phi(p)$, the inductive hypothesis gives $i_q^*\eta_p=i_q^*\gamma_p-(\gamma_p)_q u^d=0$. 
	Each class subtracted from $\gamma_p$ vanishes on every component at or above $\Phi(p)$, so $\eta_p$ has the same restrictions as $\gamma_p$ there. 
	This proves \eqref{eq:interp}, including at the other points of $U_d$ with moment value $\Phi(p)$.
\end{proof}

Normalize $\Phi$ so that $\min\Phi=0$ and let
\[
	0=t_0<t_1<\cdots<t_\ell=\max\Phi
\]
be the critical values of $\Phi$. 
Put $\Delta_i=t_{i+1}-t_i$ for $0\le i<\ell$. For each $q\in\mathcal S$, let $s(q)\in\{0,\dots,\ell\}$ be the unique index such that $\Phi(q)=t_{s(q)}$.

Now fix nonnegative integers $a,d$ satisfying \eqref{eq:range} and let $D=\max\{1,2d\}$.
For $0 \leq i < \ell$, define 
\begin{equation}\label{eq:zeta}
	\zeta_i=\sum_{\substack{p \in U_d\\\Phi(p)>t_i}}\eta_p \in H^{2d}_{S^1}(M).
\end{equation}
By \eqref{eq:interp}, the restriction of $\zeta_i$ to a point $q\in U_d$ is $u^d$ if $\Phi(q)>t_i$ and zero if $\Phi(q)\le t_i$.
We also define
\begin{equation}\label{eq:G}
	G =u^{D-2d}\sum_{i=0}^{\ell-1}\Delta_i\,\zeta_i^2-u^{D-1}\hh \in H^{2D}_{S^1}(M).
\end{equation}
\vs{0.2cm}

\begin{lemma}\label{lem:gap}
	For every $p\in U_d$, we have $G_p=0$. 
	If $q\in\mathcal S\setminus U_d$ satisfies
	\[
    		\sum_{p\in U_d}(\eta_p)_q=1,
	\]
	then $s(q)<\ell$ and $G_q\ge t_{s(q)+1}-t_{s(q)}>0$.  
\end{lemma}

\begin{proof}
	Let $p\in U_d$. By \eqref{eq:interp} and \eqref{eq:zeta}, we have
	\[
		(\zeta_i)_p = 1 \quad \text{for $t_i < t_{s(p)}$}, \qquad \text{and} \qquad (\zeta_i)_p = 0 \quad \text{for $t_i \geq t_{s(p)}$}. 
	\]
	Hence
	\[
    		G_p=\sum_{i<s(p)}\Delta_i-\Phi(p) = \sum_{i<s(p)}\Delta_i-t_{s(p)} = 0.
	\]
	This proves the first claim. 
	
	For the second claim, let $q\in\mathcal S\setminus U_d$ satisfy
	\[
		\sum_{p\in U_d}(\eta_p)_q=1. 
	\]
	If $\Phi(p)\le\Phi(q)$, or equivalently $s(p) \le s(q)$, then $(\eta_p)_q = 0$ by \eqref{eq:interp}.
	Thus some $p\in U_d$ lies strictly above $q$, and so $s(q)< s(p) \leq \ell$. 
	Moreover, for $i \leq s(q)$,
	\[
		(\zeta_i)_{q} = \sum_{\substack{p \in U_d\\ \Phi(p)>t_i}}(\eta_p)_q = \sum_{p\in U_d}(\eta_p)_q=1.
	\]
	It follows that 
	\[
    		G_q = \sum_{i=0}^{\ell-1} \Delta_i (\zeta_i)_q^2 - \Phi(q) \ge\sum_{i=0}^{s(q)}\Delta_i-t_{s(q)}=t_{s(q)+1}-t_{s(q)}>0. \qedhere 
	\]
\end{proof}

Now we are ready to prove Theorem~\ref{thm:injective}.

\begin{proof}[Proof of Theorem~\ref{thm:injective}]
Fix nonnegative integers $a,d$ satisfying \eqref{eq:range}.
Suppose $\alpha\in H^{2a}_{S^1}(M)$ is a nonzero class whose restriction to $\mathcal C_{a,d}$ vanishes.
By Kirwan's injectivity theorem, there exists an isolated fixed point $q \in \mathcal S \setminus L_d$ with $\alpha_q \neq 0$. 
Thus, the support of $\alpha$
\[
    \supp\alpha=\{q\in\mathcal S:\alpha_q\ne0\}
\]
is nonempty. 

Let $\sigma=(-1)^{n-d}$. 
At isolated points, we have $\mu(q)=n-c(q)$ so that \eqref{eq:euler} gives
\begin{equation}\label{eq:signs}
    \sgn\lambda_q=
    \begin{cases}
        \sigma,&q\in U_d\cup L_d,\\
        -\sigma,&q\in\mathcal S\setminus(U_d\cup L_d).
    \end{cases}
\end{equation}
Since $4a<2n$, applying Theorem~\ref{thm:localization} to $\alpha^2$ gives
\[
    \sum_{q\in\supp\alpha}\lambda_q\alpha_q^2=0.
\]
The sum has at least one nonzero term, so both signs appear. In particular, there exists a point  in $\supp\alpha \setminus U_d$ at which $\alpha$ is nonzero.

Set 
\[
	\beta=u^d-\sum_{p\in U_d}\eta_p \in H^{2d}_{S^1}(M).
\]
At each $q\in\mathcal S$, its coefficient is $\beta_q=1-\sum_{p\in U_d}(\eta_p)_q$. 
In particular, $\beta_q=0$ for every $q\in U_d$.
The range \eqref{eq:range} implies that $4(a+d)<2n$, so Theorem~\ref{thm:localization} applied to $\alpha^2\beta^2$ gives
\[
	0=\sum_{q\in\mathcal S}\lambda_q\alpha_q^2 \beta_q^2 =\sum_{q\in\supp\alpha\setminus U_d}\lambda_q\alpha_q^2 \beta_q^2.
\]
Every $\lambda_q$ in the last sum has sign $-\sigma$ and every $\alpha_q$ is nonzero. Hence $\beta_q=0$ at each such point, or equivalently,
\begin{equation}\label{eq:partition}
	\sum_{p\in U_d}(\eta_p)_q=1 \qquad(q\in\supp\alpha\setminus U_d).
\end{equation}

Let $G$ be the class in \eqref{eq:G}. 
By Lemma~\ref{lem:gap}, we have $G_q=0$ for every $q\in U_d$ and $G_q>0$ for $q\in\supp\alpha\setminus U_d$.
Since $4a+2D<2n$, localization for $\alpha^2G$ now gives
\[
    0=\sum_{q\in\supp\alpha\setminus U_d}\lambda_q\alpha_q^2G_q.
\]
This sum has at least one nonzero term and every term has the same sign $-\sigma$. 
Thus it gives a contradiction. 
\end{proof}

%

\begin{remark}\label{rem:positive}
For arbitrary fixed components, comparing dimensions in \eqref{eq:mixed-injection} gives
\[
	\sum_{j=0}^a b_{2j}(M)\le \sum_{F\in\cF_+}\dim H^{2a}_{S^1}(F)+|L_d|.
\]
In this bound, $\cF_+$ cannot be replaced by the components with $\mu(F)\le a$.
For example, let $S^1$ act on $\C P^4$ diagonally with weights $(0,1,2,2,3)$, and take $(a,d)=(1,0)$.
The fixed point set consists of three isolated points, corresponding to the weights $0$, $1$, $3$, and the projective line $F$ spanned by the two coordinates of weight two; note that $\mu(F)=2>a$.
The isolated points have half-coindices $4$, $3$, $0$, respectively.
Hence $L_0$ consists only of the minimum.  
Since $H^2_{S^1}(\C P^4)$ is two-dimensional, its restriction to $L_0$ cannot be injective.
Thus $F$ cannot be omitted from $\mathcal C_{1,0}$ even though $\mu(F)>a$.
\end{remark}

We now assume that all fixed points are isolated. 
To prove Theorem~\ref{thm:levels}, we first describe the classes that vanish at every fixed point with moment value at most $t$. \vs{0.2cm}

\begin{proposition}\label{prop:level-injection}
	Let $(M^{2n},\omega)$ be a closed connected Hamiltonian $S^1$-manifold with isolated fixed points and a moment map $\Phi$. 
	Given integers $a,d\ge0$ with $2a+\max\{1,2d\}<n$ and $t\in\R$, set
	\begin{equation}\label{eq:K}
		K_{a,t}=\{\alpha\in H^{2a}_{S^1}(M): \alpha_q=0\ \text{for all }q\in\mathcal S \text{ with }\Phi(q)\le t\}.
	\end{equation}
	Then the restriction map
	\[
		K_{a,t}\longrightarrow\R^{\{q\in L_d:\Phi(q)>t\}}
	\]
	is injective. Moreover, we have 
	\begin{equation}\label{eq:Kdimension}
		\dim K_{a,t}=\#\{p\in A_a:\Phi(p)>t\}.
	\end{equation}
\end{proposition}

\begin{proof}
The first claim follows directly from Theorem~\ref{thm:injective}.
For the second, recall from Lemma~\ref{lem:GT} that 
\[
	v_p=u^{a-\mu(p)}\tau_p \qquad (p\in A_a)
\]
form a basis of $H^{2a}_{S^1}(M)$. 
Each $v_p$ satisfies $i_p^*v_p\neq0$ and vanishes at every other fixed point $q$ with $\Phi(q)\le\Phi(p)$.

We claim that the $v_p$ with $\Phi(p)>t$ form a basis of $K_{a,t}$.
We first observe that they belong to $K_{a,t}$ by Lemma~\ref{lem:GT}. 
Write $\alpha\in K_{a,t}$ in the basis $(v_p)_{p\in A_a}$. 
If some $v_p$ with $\Phi(p)\le t$ has nonzero coefficient, choose such a $p$ of least moment map value. 
Since every other basis vector with nonzero coefficient vanishes at $p$, this gives $i_p^*\alpha\ne0$, contradicting $\alpha\in K_{a,t}$. This proves \eqref{eq:Kdimension}.
\end{proof}

\begin{proof}[Proof of Theorem~\ref{thm:levels}]
Inequality \eqref{eq:tail} follows by taking dimensions in Proposition~\ref{prop:level-injection}.

To prove \eqref{eq:strict-tail}, we enlarge $K_{a,t}$ by the basis vectors at level $t$ that are not indexed by $L_d$. 
Let $V\subset H^{2a}_{S^1}(M)$ be the span of the classes $v_p$ with $p\in A_a$ such that either $\Phi(p)>t$, or $\Phi(p)=t$ and $p\notin L_d$. 
Since these are distinct members of the basis in Lemma~\ref{lem:GT}, we have 
\[
	\dim V=\#\{p\in A_a:\Phi(p)>t\}+\#\{p\in A_a\setminus L_d:\Phi(p)=t\}.
\]
We claim that every class in $V$ vanishes at each $q\in L_d$ with $\Phi(q)\le t$. 
It suffices to check this for the basis vectors $v_p$ of $V$.
If $\Phi(p)>t$, then $\Phi(q)<\Phi(p)$. 
If $\Phi(p)=t$, then $q\ne p$ because $p\notin L_d$, and $\Phi(q)\le\Phi(p)$. 
In either case, $v_p$ vanishes at $q$ by Lemma~\ref{lem:GT}.

By the claim and Theorem~\ref{thm:injective}, restriction to the points of $L_d$ lying strictly above $t$ is injective on $V$. Hence
\[
	\#\{p\in A_a:\Phi(p)>t\}+\#\{p\in A_a\setminus L_d:\Phi(p)=t\}\le\#\{q\in L_d:\Phi(q)>t\}.
\]
Subtracting $\#\{q\in A_a\cap L_d:\Phi(q)>t\}$ from both sides yields
\[
	\#\{p\in A_a\setminus L_d:\Phi(p)\ge t\}\le\#\{q\in L_d\setminus A_a:\Phi(q)>t\}. 
\]
This completes the proof.
\end{proof}

The following example illustrates circle actions on complex projective space which attain equality in $|A_a|\le |L_d|$.

\begin{example}\label{ex:CPn}
Let $S^1=\R/\Z$ act on $\C P^n$ by
\[
	\theta\cdot[z_0:\cdots:z_n] =[z_0:e^{2\pi i\theta}z_1:\cdots:e^{2\pi i n\theta}z_n].
\]
The fixed points are the coordinate points $p_0,\ldots,p_n$. 
At each $p_j$, the tangent weights are $k-j$ for $k\ne j$, so $\mu(p_j)=j$ and $c(p_j)=n-j$. 
For $m\ge1$, set $a=m-1$ and take $d=0$ if $n=2m$, or $d=1$ if $n=2m+1$. 
In both cases, we have $2a+\max\{1,2d\}=n-1 < n$ and
\[
	A_a=\{p_0,\ldots,p_{m-1}\},\qquad L_d=\{p_0,p_2,\ldots,p_{2m-2}\}.
\]
Both sets have $m$ elements, so the injection in \eqref{eq:mixed-injection} is an isomorphism.
\end{example}

\begin{remark}\label{rem:range}
For the actions in Example~\ref{ex:CPn}, the map in \eqref{eq:mixed-injection} can fail to be injective when $2a+\max\{1,2d\}=n$. 
In each row below the target has smaller dimension than the source.
\[
	\begin{array}{c|c|c|c}
		M & (a,d) & \dim H^{2a}_{S^1}(M) & L_d\\ \hline
		\C P^4 & (1,1) & 2 & \{p_1\}\\
		\C P^5 & (2,0) & 3 & \{p_1,p_3\}\\
		\C P^6 & (2,1) & 3 & \{p_1,p_3\}
		\end{array}
\]
For $a\ge1$ and $d\ge2$, we do not know whether \eqref{eq:range} is optimal.
\end{remark}

\section{Betti number inequalities for isolated fixed points}\label{sec:numerical}

Throughout this section, we assume that all fixed points are isolated. 
We express both sides of $|A_a|\le |L_d|$ in terms of the even Betti numbers. \vs{0.2cm}

\begin{corollary}\label{cor:general-betti}
Let $(M^{2n},\omega)$ be a closed connected Hamiltonian $S^1$-manifold with isolated fixed points. 
For integers $a,d\ge0$ with
$2a+\max\{1,2d\}<n$,
\begin{equation}\label{eq:general-betti}
    \sum_{j=0}^{a}b_{2j}(M)
    \le \sum_{\substack{d<j\le n\\j\equiv d\pmod2}}b_{2j}(M).
\end{equation}
In particular, for $n=2m$ and $0\le d\le m-1$ one has
\begin{equation}\label{eq:even-positive}
    \sum_{j=0}^{m-d-1}b_{2j}(M)
    \le \sum_{\substack{d<j\le 2m\\j\equiv d\pmod2}}b_{2j}(M).
\end{equation}
\end{corollary}

\begin{proof}
By \eqref{eq:morse-count}, $|A_a|=\sum_{j=0}^{a} b_{2j}$.
There are $b_{2(n-j)}=b_{2j}$ points of half-coindex $j$. 
Summing over $j>d$ with $j\equiv d\pmod2$ gives $|L_d|$, so \eqref{eq:general-betti} follows from Theorem~\ref{thm:injective}. 
If $n=2m$, the choice $a=m-d-1$ satisfies \eqref{eq:range} for $0\le d\le m-1$. Then \eqref{eq:even-positive} immediately follows from \eqref{eq:general-betti}.
\end{proof}

\begin{corollary}\label{cor:odd}
Let $(M^{4m+2},\omega)$ be a closed connected Hamiltonian $S^1$-manifold
with isolated fixed points, where $m\ge1$. For $1\le d\le m$,
\begin{equation}\label{eq:odd}
    \sum_{r=0}^{\lfloor d/2\rfloor} b_{2d-4r}(M)
    \le \sum_{j=m-d+1}^{m} b_{2j}(M).
\end{equation}
In particular, $b_2(M)\le b_{2m}(M)$. 
\end{corollary}

\begin{proof}
Put $N=\sum_{j=0}^{m}b_{2j}$ and fix $1 \leq d \leq m$. Since $n=2m+1$ is odd, Poincar\'{e} duality implies that 
\[
	\sum_{\substack{0\le j\le n\\ j:\ \text{odd}}} b_{2j} =
	\sum_{\substack{0\le j\le n\\ j:\ \text{even}}} b_{2j}, 
	\qquad 
	\text{and hence} \sum_{\substack{0\le j\le n\\ j\equiv d\ (\mathrm{mod} ~2)}} b_{2j} = N.
	\qquad 
\]
On the other hand, the right-hand side of \eqref{eq:general-betti} is given by 
\[
	\sum_{\substack{d<j\le n\\ j\equiv d\ (\mathrm{mod}\ 2)}} b_{2j} = N-\sum_{\substack{0\le j\le d\\ j\equiv d\ (\mathrm{mod}\ 2)}} b_{2j} .
\]
Take $a=m-d$ so that $2a+2d=2m<n$. By definition of $N$, the left-hand side of \eqref{eq:general-betti} is
\[
	\sum_{j=0}^{m-d}b_{2j} = N-\sum_{j=m-d+1}^{m}b_{2j} .
\]
Substituting both expressions into \eqref{eq:general-betti} and cancelling $N$, we obtain
\[
	\sum_{\substack{0\le j\le d\\ j\equiv d\ (\mathrm{mod}\ 2)}} b_{2j} \le \sum_{j=m-d+1}^{m}b_{2j} ,
\]
which coincides with \eqref{eq:odd}. For $d=1$, this reads $b_2(M)\le b_{2m}(M)$.
\end{proof}

\begin{corollary}[\bibsource{Cho21}{Cho}]\label{cor:cho}
Let $(M^{4m},\omega)$ be a closed connected Hamiltonian $S^1$-manifold
with isolated fixed points, where $m\ge1$. Set $k=\lfloor m/2\rfloor$.
Then
\begin{equation}\label{eq:cho}
    \sum_{r=0}^{k-1}b_{4r+2}(M)\le\sum_{r=0}^{k-1}b_{4r+4}(M).
\end{equation}
\end{corollary}

\begin{proof}
Here $n=2m$. 
Taking $d=0$ in \eqref{eq:even-positive}, that is, $a=m-1$, gives
\[
	\sum_{j=0}^{m-1} b_{2j} \;\le \sum_{\substack{0<j\le 2m\\ j\ \mathrm{even}}} b_{2j}.
\]
The left-hand side is
\begin{equation}\label{eq:L}
	\sum_{j=0}^{m-1} b_{2j} = \sum_{\substack{0\le j<m\\ j\ \mathrm{even}}} b_{2j} \;+\; \sum_{\substack{0\le j<m\\ j\ \mathrm{odd}}} b_{2j} .
\end{equation}
On the right-hand side, the substitution $j\mapsto 2m-j$ and Poincar\'e duality give
\[
	\sum_{\substack{m<j\le 2m\\ j\ \mathrm{even}}} b_{2j}=\sum_{\substack{0\le j<m\\ j\ \mathrm{even}}} b_{2j} ,
\]
which is the even-index sum in \eqref{eq:L}.
Cancelling it, we obtain
\[
	\sum_{\substack{0\le j<m\\ j\ \mathrm{odd}}} b_{2j} \;\le \sum_{\substack{0<j\le m\\ j\ \mathrm{even}}} b_{2j} .
\]
Writing $j=2r+1$ on the left and $j=2r+2$ on the right, both sums run over $0\le r\le k-1$, where $k=\lfloor m/2\rfloor$. 
This coincides with \eqref{eq:cho}.
\end{proof}
The following table lists, for each dimension, the consequences of \eqref{eq:even-positive} and \eqref{eq:odd} with $d\ge1$ that do not follow from the others. 
For instance, in real dimension sixteen, $(a,d)=(2,1)$ gives $1+b_2+b_4\le b_6+b_{10}+b_{14}=2b_6+b_2$.
\begin{center}
\begin{tabular}{cl}
\toprule
$\dim_{\R}M$ & Consequences \\
\midrule
10 & $b_2\le b_4$ \\
14 & $b_2\le b_6$ \\
16 & $1+b_4\le2b_6$ \\
18 & $b_2\le b_8$,\quad $1+b_4\le b_6+b_8$ \\
20 & $1+b_4\le b_6+b_{10}$,\quad $b_2\le2b_8$ \\
22 & $b_2\le b_{10}$,\quad $1+b_4\le b_8+b_{10}$ \\
24 & $1+b_4+b_8\le b_6+2b_{10}$,\quad $b_2+b_6\le2b_8+b_{12}$ \\
\bottomrule
\end{tabular}
\end{center}
\vs{0.2cm}

The case $d\ge1$ also yields inequalities beyond the bounds recalled above.
For example, in real dimension sixteen, the pair $(a,d)=(2,1)$ gives $1+b_4\le 2b_6$. 
In this dimension, consider the bounds $b_2+b_6\le b_4+b_8$ \cite[Theorem~2]{Cho21}, $1+b_4\le b_2+b_6$ \cite[Corollary~7.6]{LiRen26}, and the signature bound $|b_0-b_2+b_4-\cdots+b_{16}|\le b_8$, which follows from the signature formula of Jones and Rawnsley \cite{JR}.
All three bounds hold for the formal symmetric sequence
\[
	(b_0,b_2,\ldots,b_{16})=(1,2,2,1,1,1,2,2,1),
\]
whereas $1+b_4=3>2=2b_6$. 
Hence $1+b_4\le 2b_6$ does not follow from these three bounds.

We close with real dimension twelve, where unimodality amounts to $b_2\le b_4\le b_6$. 
The first inequality is \eqref{eq:cho}, while the second lies just beyond the range \eqref{eq:range}. 
For $(a,d)=(2,1)$ we have $2a+2d=n$, and injectivity of restriction to $L_1$ would give $1+b_4\le b_6$, which fails for $\C P^6$ (Remark~\ref{rem:range}).
If the restriction of $H^4_{S^1}(M)$ to $L_1$ together with one further fixed point were injective, we would obtain $b_4\le b_6$ and hence unimodality in real dimension twelve.
For $\C P^6$, adding the minimum suffices, but we do not know whether this holds in general.

\end{document}